\documentclass[11pt]{article}
\usepackage[latin1]{inputenc}
\usepackage[english]{babel}
\usepackage{amssymb,amsfonts, amsmath, color}
\usepackage[margin=2.7cm]{geometry}
\usepackage{xcolor}
\usepackage{graphicx, color, enumerate}
\usepackage[latin1]{inputenc}
\usepackage[active]{srcltx}
\usepackage{tikz}
\usepackage{pgf}
\usepackage{etex}
\usepackage{verbatim}
\usepackage{tikz-3dplot}
\usepackage{pgfkeys}
\usepackage{fp}

\newtheorem{theorem}{Theorem}[section]
\newtheorem{proposition}[theorem]{Proposition}

\newtheorem{remark}[theorem]{Remark}
\newtheorem{definition}[theorem]{Definition}

\newcommand{\bcl}{\begin{center}}
\newcommand{\ecl}{\end{center}}
\newcommand{\brl}{\begin{right}}
\newcommand{\erl}{\end{right}}
\newcommand{\ben}{\begin{enumerate}}
\newcommand{\een}{\end{enumerate}}
\newcommand{\overliner}{\begin{array}}
\newcommand{\earr}{\end{array}}
\newcommand{\btab}{\begin{tabular}}
\newcommand{\etab}{\end{tabular}}
\newcommand{\bdoc}{\begin{document}}
\newcommand{\edoc}{\end{document}}
\newcommand{\beqy}{\begin{eqnarray}}
\newcommand{\eeqy}{\end{eqnarray}}

\newcommand{\beqi}{\begin{eqnarray*}}
\newcommand{\eeqi}{\end{eqnarray*}}
\newcommand{\bitem}{\begin{itemize}}

\newcommand{\eitem}{\end{itemize}}
\newcommand{\nln}{\newline}
\newcommand{\newt}{\newtheorem}

\newcommand{\pa}{\partial}
\newcommand{\re}{{I\!\!R}}
\newcommand{\ren}{\re^N}
\newcommand{\xr}{x\in\re }
\newcommand{\x}{\times}
\newcommand{\dyle}{\displaystyle}
\newcommand{\ene}{{I\!\!N}}
\newcommand{\irn}{\int\limits_{\re^N}}
\newcommand{\io}{\int\limits_{\O}}
\newcommand{\meas}{{\rm meas\,}}

\newcommand{\sign}{{\rm sign}}
\newcommand{\map}{\longrightarrow }
\newcommand{\imp}{\Longrightarrow }
\renewcommand{\div}{\nabla\cdot }
\newcommand{\sen}{{\rm sen\,}}
\newcommand{\tg}{{\rm tg\,}}
\newcommand{\arcsen}{{\rm arcsen\,}}
\newcommand{\arctg}{{\rm arctg\,}}
\newcommand{\supp}{{\textsl supp\ }}
\newcommand{\ity}{\int_{-\iy}^{+\iy}}
\newcommand{\limit}{\lim\limits}
\newcommand{\limi}{\limit_{n\to\infty}}
\newcommand{\sumi}{\sum\limits_{n=1}^{\infty}}
\newcommand{\ulu}{\underline u}
\newcommand{\ulw}{\underline w}
\newcommand{\ulz}{\underline z}
\newcommand{\ulv}{\underline v}
\newcommand{\uls}{\underline s}
\newcommand{\olu}{\overline u}
\newcommand{\olv}{\overline v}
\newcommand{\ols}{\overline s}
\newcommand{\ob}{\overline\b}
\newcommand{\ovar}{\overline\var}
\newcommand{\wv}{\widetilde v}
\newcommand{\wu}{\widetilde u}
\newcommand{\ws}{\widetilde s}
\renewcommand{\a }{\alpha }
\renewcommand{\b }{\beta }
\newcommand{\g }{\gamma}
\newcommand{\G }{\Gamma }
\renewcommand{\d }{\delta }

\newcommand{\D }{\Delta }
\newcommand{\e }{\varepsilon }
\newcommand{\z }{\zeta }
\renewcommand{\l }{\lambda }
\renewcommand{\L }{\Lambda }
\newcommand{\m }{\mu }
\newcommand{\n }{\nabla }
\newcommand{\s }{\sigma }
\newcommand{\Sig }{\Sigma }
\renewcommand{\t }{\tau }
\newcommand{\var }{\varphi }
\renewcommand{\o }{\omega }
\renewcommand{\O }{\Omega }
\newcommand{\bR}{{\bf R}}
\newcommand{\bC}{{\bf C}}
\newcommand{\bZ}{{\bf Z}}
\newcommand{\bN}{{\bf N}}
\newcommand{\bQ}{{\bf Q}}
\newcommand{\bK}{{\bf K}}
\newcommand{\bI}{{\bf I}}
\newcommand{\bv}{{\bf v}}
\newcommand{\bV}{{\bf V}}
\DeclareMathOperator{\suppo}{supp} \DeclareMathOperator{\di}{div}

\def\qed{\unskip\kern 6pt \penalty 500
\raise -2pt\hbox{\vrule \vbox to10pt{\hrule width 4pt
\vfill\hrule}\vrule}\par}

\newenvironment{Proof}{\removelastskip\vskip12pt
plus 1pt \noindent\em\rm}{\hfill {\qed \hskip .2cm}}

\title{Prescribed conditions at infinity \\
for fractional parabolic and elliptic equations \\ with unbounded
coefficients}
\author{Fabio
Punzo\thanks{Dipartimento di Matematica ``Federigo Enriques",
Universit\`a degli Studi di Milano, via Cesare Saldini 50, 20133
Milano, Italy (fabio.punzo@unimi.it).}\; and Enrico Valdinoci
\thanks{Weierstra{\ss} Institut f\"ur Angewandte Analysis und Stochastik,
Mohrenstra{\ss}e 39, 10117 Berlin, Germany, and Dipartimento di
Matematica ``Federigo Enriques", Universit\`a degli Studi di Milano, Via
Cesare Saldini 50, 20133 Milano, Italy, and
Istituto di Matematica Applicata e Tecnologie Informatiche ``Enrico Magenes",
Consiglio Nazionale delle Ricerche,
Via Ferrata 1, 27100 Pavia, Italy (enrico@math.utexas.edu).
Supported by the ERC grant $\epsilon$ ({\it Elliptic Pde's and
Symmetry of Interfaces and Layers for Odd Nonlinearities}). Both
authors are supported by PRIN grant 
201274FYK7 (Critical Point Theory
and Perturbative Methods for Nonlinear Differential Equations). }}

\date{}

\begin{document}
\maketitle

\abstract{We investigate existence and uniqueness of solutions to
a class of fractional parabolic equations satisfying prescribed
pointwise conditions at infinity (in space), which can be
time-dependent. Moreover, we study the asymptotic behaviour of such
solutions. We also consider solutions of elliptic
equations satisfying appropriate conditions at infinity.}


\bigskip
\smallskip

\section{Introduction}\setcounter{equation}{0}
We are concerned with existence and uniqueness of solutions to the
following linear {\it nonlocal} parabolic Cauchy problem:
\begin{equation}\label{e1}
\left\{
\begin{array}{ll}
\,   \pa_t u = -a\,(-\Delta)^{s}u + c u + f
&\textrm{in}\,\,\re^N\times (0, T]=:S_T
\\& \\
\textrm{ }u \, = u_0& \textrm{in\ \ } \re^N\times \{0\} \,,
\end{array}
\right.
\end{equation}
where the coefficient $a$ is a positive function only depending on
the space variable $x$, which becomes unbounded as $|x|\to\infty$;
$(-\Delta)^s$ denotes the fractional Laplace operator of order
$s\in (0,1), N>2s$, while $c, f, u_0 \in L^\infty(\re^N)$.
Moreover, we investigate existence and uniqueness of solutions to
the linear {\it nonlocal} elliptic equation
\begin{equation}\label{e3}
a (-\Delta)^s u - c u = f \quad \textrm{in}\;\; \re^N\,\,;
\end{equation}
in this case we also suppose that $c\leq 0$\,.

\medskip

\noindent {\bf (a) Parabolic problems\,.} The well-posedness of problem \eqref{e1} has been largely studied
in the literature in the local
case $s=1$ (see, e.g., \cite{AB}, \cite{EKP}, \cite{GMP},
\cite{IKO}, \cite{KPT}, \cite{KPu2}, \cite{KPu},
 \cite{Pinch}, \cite{Ti})\,.
As a matter of fact, if $N=1,2$ and~$s=1$, then there
exists a unique bounded solution of problem \eqref{e1}. If $N\geq
3$, a special role is played by the behaviour at infinity of the
coefficient $a$. In particular, if
\[a(x) \leq C(1+|x|^2)^{\frac{\alpha}2}\quad \textrm{for all}\;\; x\in \re^N\,,\, \textrm{for some}\,\, C>0, \a\leq 2\,,\]
then problem \eqref{e1} admits only one bounded solution (see \cite{AB}, \cite{IKO}). Instead,
if
\[
a(x) \geq C(1+|x|^2)^{\frac{\alpha}2}\quad \textrm{for all}\;\; x\in \re^N\,,\, \textrm{for some}\,\, C>0, \a>2\,,
\]
then problem \eqref{e1} admits infinitely many bounded solutions.
More precisely, for any given $g\in C([0, T])$, if
\begin{equation}\label{e8} \lim_{|x|\to \infty} u_0(x)\,=\,
g(0)\,,
\end{equation}
then there exists a unique bounded solution of problem \eqref{e1}
such that
\begin{equation}\label{e9}
\lim_{|x|\to \infty} u(x,t)=\, g(t)\, \quad \textrm{uniformly with
respect to}\,\, t\in [0, T]\,
\end{equation}
(see \cite{GMP}, \cite{KPu})\,. Observe that condition \eqref{e9} can be regarded as a Dirichlet
condition at infinity, which is time-dependent.

\medskip

More recently, existence and uniqueness results for
{\it nonlocal} Cauchy parabolic problems have been established. In this respect, in \cite{AK}, \cite{MP1},
\cite{MP2} a quite general class of integro-differential equations have been
treated; it also includes problem \eqref{e1} if there exist two constants $C_1>0, C_2>0$ such that
\begin{equation}\label{e2i}
C_1 \leq a(x) \leq C_2\quad \textrm{for all}\;\; x\in\re^N\,.
\end{equation} Furthermore,
the well-posedness for the Cauchy problem associated to the fractional porous
medium equation with a variable density $a=a(x)$ has been studied in \cite{GMP2}, \cite{GMP3},
\cite{PT}, allowing that $a(x)\to \infty$ or $a(x)\to 0$ as
$|x|\to \infty$\,. Moreover, in \cite{PVal} the uniqueness of solutions of problem \eqref{e1} with $c\equiv 0$ in suitable weighted Lebesgue spaces is stated. To be more specific,
let $$\psi(x):=(1+|x|^2)^{-\frac{\beta}2} \quad (x\in \re^N)\,,$$
$\beta$ being a positive parameter.  Suppose that, for some $C>0$ and $\a\in \re$,
$$a(x) \leq C(1+|x|^2)^{\frac{\a}2} \quad (x\in \re)\,.$$
Let $p\geq 1$, then
problem \eqref{e1} admits at most one solution $u\in L^p_\psi(S_T)$, provided that one of the next condition holds:
\begin{equation}\label{e3i}
0<\beta\leq N-2s, \a\in \re\,;
\end{equation}
\begin{equation}\label{e4i}
N-2s<\beta< N, \a\leq 2s\,;
\end{equation}
\begin{equation}\label{e5i}
\b=N, \a<2s\,;
\end{equation}
\begin{equation}\label{e6i}
\b>N, \a<2s, \a+\b<2s +N\,;
\end{equation}
here $L^p_\psi(S_T):=\big\{f: S_T\to \re\,\, \textrm{measurable such that }\,\,\int_{S_T}|f(x,t)|^p \psi(x)\, dx dt<\infty\, \big\}\,.
$
As a consequence, if $\a<2s$, we have uniqueness of solutions in the class of solutions that satisfy
\[|u(x,t)|\, \leq\, \bar C (1+|x|^2)^{\frac{\s}{2}}\quad \textrm{for all}\;\, x\in \re^N, t>0\,,\]
for some $\bar C>0$ and  $\sigma\in (0, 2s-\a)\,.$

\medskip
In the present paper, where we use completely different methods from those in \cite{PVal}, we always assume that
\[\textrm{there exist}\,\, C_0>0,\, \a>2s\,\,\textrm{such that}\;\; a(x)\geq C_0 (1+|x|^2)^{\frac{\a}2} \quad \textrm{for all}\;\; x\in \re^N\,. \leqno (H_0) \]
Clearly, this case is not covered by \cite{AK}, \cite{MP1}, \cite{MP2}, since \eqref{e2i} is not satisfied. Moreover, hypothesis $(H_0)$ excludes that \eqref{e4i}, \eqref{e5i}, \eqref{e6i} can hold; in the sequel we also discuss the case in which both $(H_0)$ and \eqref{e3i} hold.

It is worth mentioning that the unbounded diffusion coefficient $a(x)$  is very important for the applications, see for instance, for the local case, \cite{AB}, \cite{EKP}, \cite{Grig}, \cite{Pinch}, \cite{PuTe}. Clearly, the same models with the unbounded diffusion coefficient $a(x)$ occurs when considering nonlocal diffusion, for instance, in association with non-Gaussian stochastic processes, that, starting from any point in $\re^N$, can reach {\it infinity} (see, e.g., \cite{Bass}).

\smallskip

We prove (see Theorem \ref{thm1}) that there exists a unique solution of problem \eqref{e1} such that \eqref{e9} is satisfied, provided \eqref{e8} holds; furthermore,
\begin{equation}\label{e108}|u|\leq C e^{\b T}\quad \textrm{in}\;\; \re^N\times
[0,T]\,,
\end{equation}
for some $C>0$ and $\b>0$\,. This result generalizes to the case of nonlocal operator the results in \cite{GMP} and in \cite{KPu}.

In proving this result, at first for any $j\in \ene$, we consider
the {\it viscosity} solution of a suitable approximate problem in a large
cylinder~$B_j\times (0, T]$; here and hereafter for each $R>0$,
$B_R:=\big\{x\in \re^N\,:\, |x|<R \big\}$. For such problem
existence, uniqueness and regularity results have been given in
\cite{BCI1}, \cite{BCI2}. Then using suitable super- and
subsolutions and standard compactness arguments we obtain the
existence of a solution of problem \eqref{e1}, satifying the
estimate \eqref{e108}, which depends on $T$. Then, in order to
show that condition \eqref{e9} holds, proper sub- and
supersolutions are introduced (see \eqref{e26} and \eqref{e26a}
below). In the construction of these sub-- and supersolutions, which also depend on the time variable $t$ , a special role is played by a supersolution
$V\in C^2(\re^N)$ of equation
\begin{equation}\label{e18a}
-a(-\Delta)^s V\,=\, - 1 \quad \textrm{in}\,\, \re^N\setminus
\overline{B}_{R_0}\,,
\end{equation}
for some $R_0>0$, such that
\begin{equation}\label{e20}
V(x)>0 \quad \textrm{for all}\;\; x\in \re^N\,,\quad
\lim_{|x|\to\infty} V(x)\,=\,0\,,
\end{equation}
which has been appropriately constructed (see Proposition \ref{prop1}).

Moreover, we show that similar results hold for problem
\begin{equation}\label{e1b}
\left\{
\begin{array}{ll}
\,   \pa_t u = -a\,(-\Delta)^{s}u + c u + f
&\textrm{in}\,\,\re^N\times (0, \infty)
\\& \\
\textrm{ }u \, = u_0& \textrm{in\ \ } \re^N\times \{0\} \,,
\end{array}
\right.
\end{equation}
provided $c\leq 0$ (see Theorem \ref{thm1a}). Note that, in this case, condition \eqref{e9} is replaced by
\begin{equation}\label{e9a}
\lim_{|x|\to \infty} u(x,t)=\, g(t)\, \quad \textrm{uniformly with
respect to}\,\, t\in [0, \infty)\,.
\end{equation}
In order to impose condition \eqref{e9a}, we need to show preliminarily that the solution satisfies the bound
\begin{equation}\label{e109}
|u|\leq C\quad \textrm{in}\;\; \re^N\times (0,
\infty)\,,
\end{equation}
which is global in time. In order to obtain this estimate, we use a positive viscosity supersolution $h\in C(\re^N)$ of equation
\begin{equation}\label{e63}
-a(-\Delta)^s h\,=\, - 1 \quad \textrm{in}\,\, \re^N\,.
\end{equation}
Note that the proof of the
existence of such a supersolution $h$ is rather technical (see Proposition \ref{prop2}); indeed, we also show that
 \begin{equation}\label{e64}
h(x)>0 \quad \textrm{for all}\;\; x\in \re^N\,,\quad
\lim_{|x|\to\infty} h(x)\,=\,0\,.
\end{equation}

\medskip

Let us describe in general terms the deep relation between our results and stochastic calculus for jump processes. In fact, equation \eqref{e63} completed with condition \eqref{e64} can be regarded as the counterpart on $\re^N$ for the operator $a(-\Delta)^s$ of the {\it first exit-time} problem in a bounded domain for $(-\Delta)^s$. Note that the first exit-time problem in $B_R$, in the case $a\equiv 1$,  has been studied in \cite{BG}, \cite{Get}. In fact, in \cite{BG} and in \cite{Get} it is outlined the connection between the so-called {\it first exit time problem}
\begin{equation}\label{e10i} \left\{
\begin{array}{ll}
\,-(-\Delta)^{s}u  \,=\, -1 &\textrm{in}\,\,B_R
\\& \\
\textrm{ }u \, = 0 & \textrm{in\ \ }
\big(\re^N\setminus B_R\big)\,,
\end{array}
\right.
\end{equation}
and the first exit-time from $B_R$ of the jump process associated
to $(-\Delta)^s$, starting from any point in $B_R$. Moreover, it is well-known that if any point of the boundary of a
bounded domain of $\re^N$ can be reached by the jump process
associated to a nonlocal diffusion operator starting from points
inside the domain, then the Dirichlet problem admits a unique
solution that takes continuously a given datum at the boundary (see, e.g., \cite{Silv}).

Now, equation \eqref{e63}, completed with condition \eqref{e64},
corresponds to problem \eqref{e10i} in the limit case $R=\infty$, and
it is somehow related to reachability of {\it infinity} by the
jump process associated to the operator $a( -\Delta)^s$ (see
\cite{Bass}, \cite{Grig})\,. In particular, from the existence of
the supersolution $h$ it follows that {\it infinity} can actually
be attained by the jump process starting from any point $x_0\in
\re^N$. This property is usually expressed saying that the process
is {\it transient}. Therefore one can expect that
there exists a unique solution of problem \eqref{e1} which
satisfies conditions of Dirichlet type at infinity\,. Indeed, we
prove this.

We should mention that, to the best of our knowledge, in the
literature no results concerning the prescription of general
Dirichlet conditions at infinity for solutions of nonlocal
parabolic (or elliptic) equations have been obtained before the present paper\,.
\medskip

Finally, we prove that the solution $u(x,t)$ of problem \eqref{e1b} satisfying \eqref{e9a} admits a limit function as $t\to \infty$. In fact, the function
\[W(x):= \lim_{t\to \infty} u(x,t) \quad (x\in \re^N)\]
is the unique solution of equation \eqref{e3} such that
\begin{equation}\label{e1f}
\lim_{|x|\to \infty} W(x)\,=\, \g\,,
\end{equation}
provided
\begin{equation}\label{e11}
\g=\lim_{t\to \infty} g(t)
\end{equation}
(see Theorem \ref{thm3}). Such result is shown by adapting to the present situation the method of sub- and supersolutions used in \cite{Satt} in the case of bounded domains of $\re^N$ for "local" parabolic equations. Indeed, some important changes are in order, in view of the nonlocal character of the problem and since we prescribe conditions as $|x|\to \infty$\,.

\smallskip

Now, let us discuss the case that both $(H_0)$ and \eqref{e3i} hold. In view of existence results described above, for any $g_1, g_2\in C([0, T])$ with $g_1\not \equiv g_2$ there exist two solutions $u_1$ and $u_2$ of problem \eqref{e1} such that
\[ u_1(x, t)\to g_1(t) \,, \quad u_2(x,t)\to g_2(t) \quad \textrm{as}\;\; |x|\to \infty, \quad \textrm{uniformly for}\,\, t\in [0, T]\,.\]
Set $w:= u_1-u_2$. Since $g_1\not \equiv g_2$, there exists $t_0\in [0, T]$ such that
\[ w(x, t_0) \to g_1(t_0)-g_2(t_0)\neq 0\quad \textrm{as}\;\; |x|\to \infty\,.\]
Therefore, $w\not \in L^p_\psi(S_T)$, with the choice of $\b$ required in \eqref{e3i}\,. Hence, the uniqueness result in \cite{PVal} cannot be applied to conclude that $w\equiv 0$. So, obviously, the results in \cite{PVal} and those described above are not in contradiction.
\bigskip

\noindent {\bf (b) Ellipitc equations.} In the local case,
some existence and uniqueness results for equations \eqref{e3} with $s=1$ can be deduced from general results in \cite{PuTe}\,. Moreover, the case $0<s<1$ has been treated in \cite{PVal}\,. In particular, it is shown that uniqueness results in $L^p_\psi(\re^N)$, analogous to those recalled above for the parabolic problem, holds, if $c\leq - c_0$ with $c_0>0$ large enough\,. Consequently, if $\a<2s$, we have uniqueness of solutions in the class of solutions that satisfy
\[|u(x)|\, \leq\, \bar C (1+|x|^2)^{\frac{\s}{2}}\quad \textrm{for all}\;\, x\in \re^N\,,\]
for some $\bar C>0$ and  $\sigma\in (0, 2s-\a)\,.$ On the other hand, only requiring that $c_0>0$, it is shown uniqueness in $L^p_{(1+|x|)^{N-2s+\a}}(\re^N)$, if $\a<2s$\,.

From the result concerning the asymptotic behaviour of solutions of problem \eqref{e1b} recalled in $(a)$ above, we can infer that for any $\g\in \re$ there exists a unique solution $u$ of equation \eqref{e3}, which satisfies 
\begin{equation}\label{e10c}
\lim_{|x|\to \infty} u(x)\,=\, \g\, \,.
\end{equation}
However, we also prove this existence and uniqueness result independently, without using results for parabolic problems.
In fact, we solve a proper approximate problem in a large ball~$B_j$ for any $j\in \ene$. In order to obtain a uniform bound, for any $j\in \mathbb N$, for the solutions of such problems we use in crucial way the supersolution $h$ of equation \eqref{e63}. Then, by standard compactness tools, we get a solution of equation \eqref{e3}. Using again the supersolution $h$, and in particular the fact that \eqref{e64} holds, we impose that
\eqref{e10c} holds.

We devote the forthcoming Section~\ref{mf} to the precise
statement of the main results obtained in this paper
(see in particular Subsection~\ref{MS:RES:90}).

\section{Mathematical framework and results}\label{mf}\setcounter{equation}{0} \label{mfr}
The fractional Laplacian $(-\Delta)^{s}$ can be defined by Fourier
transform. Namely, for any function $g$ in the Schwartz class
$\mathcal S$, we say that
$$(-\Delta)^{\sigma/2} g = h\,,$$
if
\begin{equation}
\label{02015} \hat{h}(\xi)= |\xi|^{\sigma}\hat{g}(\xi).
\end{equation}
Here, we used the notation~$\hat{h}={\mathfrak F} h$ for the
Fourier transform of~$h$. Furthermore, consider the space
\[\mathcal L^s(\re^N):= \left\{u:\re^N\to \re\,\;\textrm{measurable}\, \,|\,\int_{\re^N}\frac{|u(x)|}{1+|x|^{n+2s}}dx<\infty \right\}\,, \]
endowed with the norm
\[ \| u \|_{\mathcal L^s(\re^N)}:=\int_{\re^N}\frac{|u(x)|}{1+|x|^{N+2s}}dx\,.\]
If $u\in \mathcal L^s(\re^N)$ (see \cite{Silv}), then $(-\Delta)^s
u$ can be defined as a distribution, $i.e.$, for any $\varphi\in
\mathcal S$,
\[ \int_{\re^N} \varphi (-\Delta)^s u \, dx \,=\,   \int_{\re^N}  u (-\Delta)^s \varphi \, dx\,. \]
In addition, suppose that, for some $\g>0$, $u\in \mathcal
L^s(\re^N)\cap C^{2s+\g}(\re^N)$ if $s<\frac 1 2$, or  $u\in
\mathcal L^s(\re^N)\cap C^{1, 2s+\g-1}_{loc}(\re^N)\;
\textrm{if}\; s\geq \frac 1 2$. Then we have
\begin{equation}
\label{ea1} (-\Delta)^{s} u(x)=C_{N,s}\,\, \textrm{P.V.}\,
\int_{\re^N} \frac{u(x)-u(y)}{|x-y|^{N+2s}}d y\quad (x\in \re^N),
\end{equation}
where $$C_{N,s}=\frac{2^{2s}{s}
\Gamma((N+2s)/2)}{\pi^{N/2}\Gamma(1-s)},$$ $\Gamma$ being the
Gamma function; moreover, $(-\Delta)^s u\in C(\re^N)$. In the
sequel, for simplicity, we shall write \[\int_{\re^N}
\frac{u(x)-u(y)}{|x-y|^{N+2s}}d y \equiv \textrm{P.V.}\,
\int_{\re^N} \frac{u(x)-u(y)}{|x-y|^{N+2s}}d y   \quad (x\in
\re^N)\,.\]

Note that the constant $C_{N,s}$ satisfies the identity
\[ (-\Delta)^s u = \mathfrak F^{-1} \big( |\xi|^{2s}\mathfrak F u  \big)\,,\quad \xi\in\re^N, u\in \mathcal S\,,\]
so (see \cite{DPVal})
\[C_{N,s}=\left(\int_{\re^N}\frac{1-\cos(\xi_1)}{|\xi|^{N+2s}}d\xi \right)^{-1}\,.\]

\medskip

Concerning the coefficients $a$ and $c$, and the function $f$  we always make the following assumption:
\[
\textrm{\ \ } \left\{
\begin{array}{l}
(i) \quad \;  a\in C^{0,\s}_{loc}(\re^N)\,\,\big(\s\in
(0,1)\big),\;\; a(x)>0\quad \textrm{for all}\;\; x\in \re^N\,;
\\
(ii)\quad \! c, f\in C^{0,\s}_{loc}(\re^N)\cap L^\infty(\re^N)\,.
\end{array}
\right. \leqno(H_1) \]

Now we can give the definition of solution. Let $\Omega\subseteq
\re^N$ be an open subset.

\begin{definition}\label{defsoleqp}
We say that a function $u$ is a {\em subsolution} to equation
\begin{equation}\label{e3a}
\pa_t u\,= -\,a\, (-\Delta)^s u\, + c u + f\quad \textrm{in}\;\;
Q_T:=\Omega\times(0, T]\,,
\end{equation}
if
\begin{itemize}
\item[(i)] $u$ is upper semicontinuous in $S_T$\,; \item[(ii)] for
any open bounded subset $U\subset Q_T$, for any $(x_0, t_0)\in U$,
for any test function $\varphi\in C^2(S_T)$ such that $u(x_0,
t_0)-\varphi(x_0, t_0)\geq u(x,t)-\varphi(x,t)$  for all $(x,t)\in
U$, one has
\[ \pa_t \psi(x_0, t_0) \leq - a(x_0)(-\Delta)^s \psi(x_0, t_0) \,+ c(x_0) u(x_0, t_0) + f(x_0)\,,\]
where \begin{equation}\label{e10}
\psi:=\left\{
\begin{array}{ll}
\,\varphi &  \textrm{in}\,\,U
\\& \\
\textrm{ }u & \textrm{in\ \ } S_T\setminus U \,.
\end{array}
\right.
\end{equation}
\end{itemize}

Furthermore, we say that a function $u$ is a {\em supersolution}
to equation \eqref{e3a} if
\begin{itemize}
\item[(i)] $u$ is lower semicontinuous in $S_T$\,; \item[(ii)] for
any open bounded subset $U\subset Q_T$, for any $(x_0, t_0)\in U$,
for any test function $\varphi\in C^2(S_T)$ such that $u(x_0,
t_0)-\varphi(x_0, t_0)\leq u(x,t)-\varphi(x,t)$ for all $(x,t)\in
U$, one has
\[ \pa_t \psi(x_0, t_0) \geq - a(x_0)(-\Delta)^s \psi(x_0, t_0)\,+ c(x_0) u(x_0, t_0) + f(x_0)\,,\]
\end{itemize}
where $\psi$ is defined by \eqref{e10}. Finally, we say that $u$
is a {\em solution} to equation \eqref{e3} if it is both a
subsolution and a supersolution to equation \eqref{e3a}.
\end{definition}

Let $g\in C([0,T]),\, u_0\in C(\re^N)$ with
\begin{equation}\label{e7}
u_0(x, 0)=g(0)\quad \textrm{for all} \;\;\, x\in \re^N\setminus
\Omega.
\end{equation}
Consider the problem
\begin{equation}\label{e4} \left\{
\begin{array}{ll}
\,\pa_t u = -a\, (-\Delta)^{s}u  + c u + f &\textrm{in}\,\,Q_T
\\& \\
\textrm{ }u \, = g & \textrm{in\ \ }
\big(\re^N\setminus\Omega\big)\times (0, T] \,
\\& \\
\textrm{ }u \, = u_0& \textrm{in\ \ } \re^N\times \{0\} \,.
\end{array}
\right.
\end{equation}

\begin{definition}\label{defsolp}
We say that a function $u$ is a {\em subsolution} to problem
\eqref{e4} if
\begin{itemize}
\item[(i)] $u$ is upper semicontinuous in $\overline{S_T}$\,;
\item[(ii)] $u$ is a subsolution to equation \eqref{e3a}\,;
\item[(iii)] $u(x, t)\leq g (t)$ for all $x\in
\re^N\setminus\Omega, t\in (0, T]$ and $u(x,0)\leq u_0(x)$\;\; for
all\;\; $x\in \re^N\,.$
\end{itemize}
Similarly, supersolutions are defined. Finally, we say that $u$ is
a {\em solution} to problem \eqref{e4} if it is both a subsolution
and a supersolution to problem \eqref{e4}.
\end{definition}
Observe that according to our definition, any solution of problem
\eqref{e4} takes continuously the initial datum $u_0$ and the
boundary datum $g$.

\begin{definition}\label{defsole}
We say that a function $u$ is a {\em subsolution} to equation
\begin{equation}\label{e5}
a\,(-\Delta)^s u\, - c u = f\quad \textrm{in}\;\; \Omega\,,
\end{equation}
if
\begin{itemize}
\item[(i)] $u$ is upper semicontinuous in $\re^N$\,; \item[(ii)]
for any open bounded subset $U\subset \Omega$, for any $x_0\in U$,
for any test function $\varphi\in C^2(\re^N)$ such that
$u(x_0)-\varphi(x_0)\geq u(x)-\varphi(x)$ for all $x\in U$, one
has
\[ a(x_0)(-\Delta)^s \psi(x_0) \,- c(x_0) u(x_0)\,\leq \,f(x_0)\,,\]
\end{itemize}
where $\psi$ is defined by
\begin{equation}\label{e10a}
\psi:=\left\{
\begin{array}{ll}
\,\varphi &  \textrm{in}\,\,U
\\& \\
\textrm{ }u & \textrm{in\ \ } \re^N\setminus U \,.
\end{array}
\right.
\end{equation}

Furthermore, we say that a function $u$ is a {\em supersolution}
to equation \eqref{e5} if
\begin{itemize}
\item[(i)] $u$ is lower semicontinuous in $\re^N$\,; \item[(ii)]
for any open subset $U\in \Omega$, for any $x_0\in U$, for any
test function $\varphi\in C^2(\Omega)$ such that
$u(x_0)-\varphi(x_0)\leq u(x)- \varphi(x)$ for all $x\in U$, one
has one has
\[   a(x_0)(-\Delta)^s \psi(x_0)\,- c(x_0) u(x_0)\, \geq \, f(x_0)\,.\]
\end{itemize}

Finally, we say that $u$ is a {\em solution} to equation
\eqref{e5} if it is both a subsolution and a supersolution to
equation \eqref{e5}.
\end{definition}

Consider the following problem
\begin{equation}\label{e6} \left\{
\begin{array}{ll}
\, -a\, (-\Delta)^{s}u  - c u \,=\, f &\textrm{in}\,\,\Omega
\\& \\
\textrm{ }u \, = \g & \textrm{in\ \ }
\big(\re^N\setminus\Omega\big)\,,
\end{array}
\right.
\end{equation}
where $\g\in \re$.

\begin{definition}\label{defsolpe}
We say that a function $u$ is a {\em subsolution} to problem
\eqref{e6} if
\begin{itemize}
\item[(i)] $u$ is upper semicontinuous in $\re^N$\,; \item[(ii)]
$u$ is a subsolution to equation \eqref{e5}\,; \item[(iii)]
$u(x)\leq \g$ for all $x\in \re^N\setminus\Omega.$
\end{itemize}
Similarly, supersolutions and solutions are defined.
\end{definition}

In the next two Remarks we summarize existence, uniqueness
and regularity results shown in \cite{BCI1, BCI2}, for problems
\eqref{e4} and \eqref{e6}, that will be used in the sequel.

\begin{remark}\label{oss1} Let $\Omega\subset\re^N$ be an open bounded subset with $\pa
\Omega$ of class $C^1$. Let assumption $(H_1)$ be satisfied. Let
$g\in C([0,T]), u_0\in C(\re^N)\cap L^\infty(\re^N)$; suppose that
condition \eqref{e7} is satisfied. We have that
\begin{itemize}
\item[$(i)$] there exists a unique solution to problem \eqref{e4};
\item[$(ii)$] the comparison principle holds for problem
\eqref{e4}; \item[$(iii)$] if $u$ is a solution of equation
\eqref{e3a}, then, for some $0<\mu<1,$ for any open subset
$\Omega'\subset\subset\Omega, \tau\in (0, T]$ we have
\[|u(x, t_1)-u(y, t_2)|\leq C\big(|x-y|^\mu + |t_1-t_2|^{\frac{\mu}{2s}}\big)\quad \textrm{for all}\;\; x,y\in \Omega', t_1,t_2\in [\t, T],\]
for some constant $C>0$, which only depends on $\|u\|_\infty, N,
a, c, f$\,.
\end{itemize}
Note that $(i)-(ii)$ follow from \cite[Section 4.3]{BCI1}, while
$(iii)$ is a consequence of the results in \cite{BCI2} (see also the comments at the
end of page 2 in \cite{BCI2}).
\end{remark}

\begin{remark}\label{oss2}  Let $\Omega\subset\re^N$ an open bounded subset with $\pa
\Omega$ of class $C^1$; let $\g\in \re$. Let assumption $(H_1)$ be
satisfied. Assume that $c\leq 0$ in $\Omega$. We have that
\begin{itemize}
\item[$(i)$] there exists a unique solution to problem \eqref{e6};
\item[$(ii)$] the comparison principle holds for problem
\eqref{e6}; \item[$(iii)$] if $u$ is a solution of equation
\eqref{e5}, then, for some $\mu\in (0,1),$ for any open subset
$\Omega'\subset\subset\Omega$,
\[\| u\|_{C^{0, \mu}(\Omega')}\leq C\,,\]
for some constant $C>0$, which only depends on $\|u\|_\infty, N,
a, c, f$.
\end{itemize}
Note that $(i), (ii)$ follow from \cite[Theorem 2]{BCI1}) and
\cite[Theorem1]{BCI2}), whereas from Theorem 2 and the comments at
the end of page 2 in \cite{BCI2} it follows $(iii)$.
\end{remark}

\subsection{Main  results: existence, uniqueness and asymptotic behaviour of solutions}\label{MS:RES:90}
Concerning existence and uniqueness of solutions of problem \eqref{e1} we have the next result.
\begin{theorem}\label{thm1}
Let assumptions $(H_0), (H_1)$ be satisfied. Let $T>0$. Let $g\in C([0,T]),
u_0\in C(\re^N)\cap L^\infty(\re^N)$; suppose that condition
\eqref{e8} is satisfied. Then there exists a unique solution $u$
to problem \eqref{e1} such that condition \eqref{e9} is satisfied.
Furthermore, \eqref{e108} holds\,.
\end{theorem}

Under the extra hypothesis that $c\leq 0$, we have the next existence and uniqueness for problem \eqref{e1b}\,.
\begin{theorem}\label{thm1a}
Let assumptions $(H_0), (H_1)$ be satisfied. Let $g\in
C([0,\infty))\cap L^\infty((0,\infty)), u_0\in C(\re^N)\cap
L^\infty(\re^N), c\leq 0$; suppose that condition \eqref{e8} is
satisfied. Then there exists a unique solution to problem
\eqref{e1b} such that condition \eqref{e9a} is satisfied.
Furthermore, for some $C>0$, \eqref{e109} holds\,.
\end{theorem}

\begin{remark}
Observe that the estimate in \eqref{e108} depends on $T>0$, while
that in \eqref{e109} is independent of $T$. In order to get
\eqref{e109} we use the further hypothesis $c\leq 0$.
\end{remark}

Concerning the elliptic equation \eqref{e3} we show the next result.
\begin{theorem}\label{thm2}
Let assumptions $(H_0), (H_1)$ be satisfied. Let $\g\in \re$;
suppose that $c\leq 0$ in $\re^N$. Then there exists a unique
solution to equation \eqref{e3} such that condition \eqref{e10c}
is satisfied.
\end{theorem}

The next theorem is concerned with the asymptotic behaviour as $t\to \infty$ of solutions of problem \eqref{e1b}\,.
\begin{theorem}\label{thm3}
Let assumptions of Theorem \ref{thm1a} be satisfied. Let
$\g:=\lim_{t\to \infty} g(t)\,.$ Let $u$ be the unique solution to
problem \eqref{e1} such that \eqref{e9} is satisfied. Suppose that
condition \eqref{e11} holds. Then
\[\lim_{t\to \infty} u(x,t)\,=\, W(x)\quad \textrm{for all}\;\; x\in \re^N\,,\]
where $W$ is the unique solution of equation \eqref{e3} satisfying
condition \eqref{e1f}.
\end{theorem}

\begin{remark}
Note that the existence result in Theorem \ref{thm2} can be regarded as a consequence of
Theorem \ref{thm3}. In fact, from Theorem \ref{thm3} in particular
we obtain the existence of a solution $W(x):=\lim_{t\to \infty}
u(x,t)$ of problem \eqref{e3}, where $u(x,t)$ is the solution of
problem \eqref{e1} with $g(t)\equiv \g\,$ and $u_0$ satisfying
\eqref{e8}. However, in Section \ref{pp} we
give an independent proof of
Theorem \ref{thm3}, without using results concerning the parabolic
problem. Finally, observe that the supersolution $h(x)$ of
equation \eqref{e63} plays a crucial role both in the Proof of
Theorem \ref{thm2} and in that of Theorem \ref{thm3}\,.
\end{remark}

\section{Construction of stationary supersolutions}
Let us introduce the hypergeometric function
\[ _2 F_1 (a,b,c,\s)\equiv F(a,b,c,\s)\]
with $a,b\in \re, c>0,\s\in \re\setminus\{1\}\,$. The next limits
holds (see \cite[Chapters 15.2, 15.4]{DLMF}):
\begin{equation}\label{e16}
\lim_{\s\to 1^-}
F(a,b,c,\s)=\frac{\G(c)\G(c-a-b)}{\G(c-a)\G(c-b)}\,,
\end{equation}
where $\G$ is the Gamma function. Note that
\begin{equation}\label{e17}
\Gamma(t)>0 \quad \textrm{for all}\;\; t>0\,.
\end{equation}

For any $C>0, \b>0$ define the function
\begin{equation}\label{e19}
V(x):= C(1+|x|^2)^{-\frac{\b}2}\quad \big(x\in
\re^N\big)\,.\end{equation}

Concerning the function $V$, we show the next result.
\begin{proposition}\label{prop1}
Let assumptions $(H_0), (H_1)-(i)$ be satisfied. There exists
$C>0, \b>0$ and $R_0>0$ such that the function $V$ satisfies
\begin{equation}\label{e18}
-a(x)(-\Delta)^s V(x)\,\leq\, - 1 \quad \textrm{for all}\,\, x\in
\re^N\setminus \overline{B}_{R_0}\,.
\end{equation}
In particular, $V$ is a supersolution of equation \eqref{e18a}
in the sense of Definition \ref{defsole}\,. Moreover, \eqref{e20} holds\,.
\end{proposition}

\noindent{\it Proof\,.\,} To begin with, observe that since $V\in
C^\infty(\re^N)$, we have that $(-\Delta)^s V\in C(\re^N)$ (see
Section \ref{mf}). From the proof of Corollary 4.1 in
\cite{FerrVerb} it follows that, for some constant $\check C>0$,
\begin{equation}\label{e21}
-(-\Delta)^s V(x)\,=\,- C \check C F(a,b,c, -|x|^2)\quad
\textrm{whenever}\;\, |x|>1\,,
\end{equation}
with
\[a=\frac N2 +s,\quad b=\frac{\b}2 +s, \quad c=\frac N2\,.
\]
By Pfaff's transformation,
\begin{equation}\label{e22}
F(a,b,c, -|x|^2)\,=\; \frac 1{(1+|x|^2)^b} F\left(c-a, b, c,
\frac{|x|^2}{1+|x|^2}\right)\quad \textrm{for all}\;\; x\in
\re^N\setminus B_1\,.
\end{equation}
Suppose that $0<\b<N$. Hence, as a consequence of \eqref{e16} we
get
\begin{equation}\label{e23}
\lim_{|x|\to\infty} F\left(-s, \frac{\b}2+s, \frac N2,
\frac{|x|^2}{1+|x|^2}\right)\,=\,\frac{\Gamma\left(\frac
N2\right)\Gamma\left(\frac{N-\b}2\right)}{\Gamma\left(\frac N2
+s\right)\Gamma\left(\frac{N-\b}2-s\right)}=:K\,.
\end{equation}
Now, we choose $0<\b<N-2s$, so we have $K>0$. Due to \eqref{e21},
\eqref{e22}, \eqref{e23}, we can find $R_0>0$ such that
\begin{equation}\label{e24}
-(-\Delta)^s V(x) \leq -\frac{C \check C K }{2(1+|x|^2)^{\frac
{\b}2+s}}\quad \textrm{for all}\;\; x\in \re^N\setminus B_{R_0}\,.
\end{equation}
If we select $\b$ that also satisfies $0<\b\leq \a-2s$, then
\eqref{e24} and $(H_0)$ yields \eqref{e18}, provided that
\[C\geq \frac 2{ C_0 \check C K}\,.\]

Since $V\in C^\infty(\re^N)$, it easily follows that $V$ is a
supersolution of equation \eqref{e18} in the sense of Definition
\ref{defsole}. Finally, the properties in \eqref{e20} immediately
follow from the very definition of $V$. \hfill $\square$

\begin{proposition}\label{prop2}
Let assumptions $(H_0), (H_1)-(i)$ be satisfied. Then there exists a
supersolution $h$ of equation \eqref{e63} in the sense of Definition \ref{defsole}, which satisfies
\eqref{e64}\,.
\end{proposition}

\noindent{\it Proof\,.\,\,} Let $R_0, C$ and $V$ be given by
Proposition \ref{prop1}. Take $\hat R>R_0$. From the results in
\cite{Get} it follows that, for a certain $C_1=C_1(N,s)>0$, the
function
\[\hat W(x)\equiv \hat W(|x|):=  C_1\big(\hat R^2 - |x|^2 \big)_+^{s/2}\quad (x\in \re^N)\,.\]
solves
\begin{equation}\label{e65} \left\{
\begin{array}{ll}
-(-\Delta)^{s}u \, =\, - 1 &\textrm{in}\,\, B_{\hat R}
\\& \\
\textrm{ }u \, = 0 & \textrm{in\ \ } \re^N\setminus B_{\hat
R}\,\,.
\end{array}
\right.
\end{equation}
Hence, it easily follows that for each $\mu_0>0, \mu_1\geq \mu_0
\max_{\overline B_{\hat R}}\frac 1 a$ and $\mu_2>0$, the function
\[W(x)\equiv W(|x|):= \mu_1\hat W(x) + \mu_2\quad \big(x\in \re^N\big)\]
is a supersolution of problem
\begin{equation}\label{e66} \left\{
\begin{array}{ll}
-a(x)(-\Delta)^{s}u \, =\,- \mu_0 &\textrm{in}\,\, B_{\hat R}
\\& \\
\textrm{ }u \, = \mu_2 & \textrm{in\ \ } \re^N\setminus B_{\hat
R}\,\,.
\end{array}
\right.
\end{equation}
For any $\tilde C>0$ set
\[\tilde V(x)\equiv \tilde V(|x|):= \tilde C V(|x|) \quad \big( x\in \re^N\big)\,.\]
We see that for suitable $\mu_2>0, \hat R>0, \tilde C>0$, possibly
depending on $C, C_1, \mu_1, \beta, s$, we have
\begin{equation}\label{e75b}
\tilde V(0) > W(0)\,,\quad \tilde V < W \quad \textrm{in}\;\;
\left[\frac{\hat R}2, \hat R \right]\,.
\end{equation}
In fact, if
\begin{equation}\label{e76b}
\tilde C C\left[1+\left(\frac{\hat
R}2\right)^2\right]^{-\frac{\b}2} < \mu_2 < \tilde C C- \mu_1 C_1
\hat R^s\,,
\end{equation}
then  \eqref{e75b} holds.  Now, if we take
\begin{equation}\label{e76c}\tilde C=\hat R^{s+\nu},\quad
\mu_2=\hat R^{s+\nu-\b+\d}\end{equation} with
\[0<\d<\b,\quad \nu>\max\{0, \b-s+2\},\]
then, it is direct to see that, for $\hat R>0$ large enough,
\eqref{e76b}, and so \eqref{e75b}, holds.

In view of \eqref{e75b}, there exists $\bar R\in (0, \hat R/2)$
such that $W(\bar R)=\tilde V(\bar R)\,.$ Indeed, such $\bar R$ is
unique. To see this, take any $\bar R>0$ such that $W(\bar
R)=\tilde V(\bar R)\,.$ In view of \eqref{e75b} and the very
definition of $W$ and $\tilde V$ we have that $\bar R\in (0, \hat
R/2)\,.$ So,
\begin{equation}\label{e150}
\hat R^2-\bar R^2 \geq 1\,,
\end{equation}
provided $\hat R>2\,.$ Moreover, it is direct to check that if we
show that
\begin{equation}\label{e88}
W'(\bar R) > \tilde V'(\bar R)\,,
\end{equation}
then such $\bar R$ is unique. In order to show \eqref{e88}, note
that \eqref{e88} is equivalent to
\begin{equation}\label{e89}
s\mu_1 C_1(1+\bar R^2)<\b (\hat R^2-\bar R^2)^{1-\frac s 2}\tilde
V(\bar R)\,.
\end{equation}
Now, since
\[s\mu_1 C_1(1+\bar R^2)\leq s \mu_1 C_1(1+\hat R^2)\,,\]
in view of \eqref{e150}, \eqref{e89} follows if we prove that
\begin{equation}\label{e90}s\mu_1C_1(1+\hat R^2)\leq \b \mu_2\,.\end{equation} Clearly,
\eqref{e90} is a direct consequence of \eqref{e76c}, provided that
$\hat R>0$ is large enough. Hence, we have that \eqref{e88} is
satisfied, and so $\bar R$ is unique. Therefore,
\begin{equation}\label{e67}
\tilde V\geq W \quad \textrm{in}\;\; B_{\bar R},\quad \tilde
V(\bar R)=W(\bar R), \quad \tilde V\leq W \quad \textrm{in}\;\;
\re^N\setminus B_{\bar R}\,.
\end{equation}
Furthermore, since
\[\tilde V(\bar R)=W(\bar R),\]
we get
\[\mu_2(1+\bar R)^{\frac{\b}2}\geq \tilde C C - \mu_1 C_1 \hat R^s (1+\bar R^2)^{\frac{\b}2},\]
thus, \eqref{e76c} yields
$$
(1+\bar R^2)^{\frac{\b}2} \geq \frac{C \hat R^{s+\nu}}{\hat
R^{s+\nu+\delta-\b}+\mu_1 C_1 \hat R^s}\,.
$$
This implies that we can choose $\hat R>0$ so large that $\bar
R\in (R_0, \hat R/2)\,.$

\medskip

Define
\[h:= \min\{\tilde V,\, W \}\quad \textrm{in}\;\; \re^N\,.\]
We claim that $h$ is a supersolution of equation
\[-a(-\Delta)^s h = -\min\{\mu_0, \tilde C\}\quad \textrm{in}\;\; \re^N\,.
\]
In fact, since $V$  is a supersolution of equation \eqref{e18a},
by Definition \ref{defsole} and \eqref{ea1}, for any open bounded
subset $\Omega'\subset \re^N\setminus \overline{B}_{R_0}$, for any
$x_0\in \Omega'$, for any test function $\varphi\in C^2(\re^N)$
such that $\tilde V(x_0)-\varphi(x_0)\leq \tilde V(x)-\varphi(x)$ for
all $x\in \Omega'$, one has
\[ a(x_0)C_{N,s}\,\,
\int_{\re^N} \frac{\psi(x_0)-\psi(y)}{|x_0-y|^{N+2s}}d y \,\geq
\tilde C \,\,,\] where $\psi$ is defined by \eqref{e10a} with $u$
replaced by $\tilde V$ and $U$ by $\Omega'$. Hence
\begin{equation}\label{e68}
a(x_0) C_{N,s}\left\{\int_{\Omega'}
\frac{\varphi(x_0)-\varphi(y)}{|x_0-y|^{N+2s}}d y +
\int_{\re^N\setminus \Omega'} \frac{\tilde V(x_0)-\tilde
V(y)}{|x_0-y|^{N+2s}}d y\right\}\geq  \tilde C\,.
\end{equation}

Similarly, since $W$ is a supersolution of problem \eqref{e66}, we
have that for any open bounded subset $U\subset B_{R_0}$, for any
$x_0\in U$, for any test function $\varphi\in C^2(\re^N)$ such
that $W(x_0)-\varphi(x_0)\leq W(x)-\varphi(x)$ for all $x\in U$,
one has
\begin{equation}\label{e69}
a(x_0) C_{N,s}\left\{\int_{U}
\frac{\varphi(x_0)-\varphi(y)}{|x_0-y|^{N+2s}}d y +
\int_{\re^N\setminus U} \frac{W(x_0)- W(y)}{|x_0-y|^{N+2s}}d
y\right\}\geq  \mu_0\,.
\end{equation}

\medskip

Now, take any $x_0\in \re^N$ with $|x_0|\geq \bar R $, any open
bounded subset $U\subset \re^N$ with $x_0\in U$, and any test
function $\varphi\in C^2(\re^N)$ such that
$h(x_0)-\varphi(x_0)\leq h(x)-\varphi(x)$ for all $x\in U$. Set
\begin{equation}\label{e73} \psi:=\left\{
\begin{array}{ll}
\,\varphi &  \textrm{in}\,\,U
\\& \\
\textrm{ }h & \textrm{in\ \ } \re^N\setminus U \,.
\end{array}
\right.\end{equation}

Note that, due to \eqref{e67}, we have
\begin{equation}\label{e71}
h(x_0)= \tilde V(x_0)\,.
\end{equation}
For any $0< \epsilon< \bar R- R_0$, we have $\mathcal U_1:= U\cap
\big(\re^N\setminus B_{R_0+\epsilon}\big)\subset \re^N\setminus
\overline{B}_{R_0}, x_0\in \mathcal U_1$. Moreover,
\begin{equation}\label{e71a}\varphi(x)\leq \tilde V(x) \quad \textrm{for all}\;\; x\in
\mathcal U_1, \quad \varphi(x_0)=\tilde V(x_0)\,.\end{equation}
So, from \eqref{e68} with $\Omega'=\mathcal U_1$ we get
\begin{equation}\label{e70}
a(x_0) C_{N,s}\left\{\int_{\mathcal U_1}
\frac{\varphi(x_0)-\varphi(y)}{|x_0-y|^{N+2s}}d y +
\int_{\re^N\setminus\mathcal U_1} \frac{\tilde V(x_0)-\tilde
V(y)}{|x_0-y|^{N+2s}}d y\right\}\geq  \tilde C\,.
\end{equation}
Due to \eqref{e67} and \eqref{e71}, since $h\leq \tilde V$ in
$\re^N$, we have
\begin{equation}\label{e72}
a(x_0) C_{N,s}\left\{\int_{\mathcal U_1}
\frac{\varphi(x_0)-\varphi(y)}{|x_0-y|^{N+2s}}d y +
\int_{\re^N\setminus \mathcal U_1} \frac{h(x_0)-
h(y)}{|x_0-y|^{N+2s}}d y\right\}\geq  \tilde C\,.
\end{equation}
Set $\mathcal U_2:= U\cap B_{R_0+\epsilon}$.\, In view of
\eqref{e72}, since $\varphi(x_0)-\varphi(y)\geq h(x_0)- h(y)$ for
all $y\in \mathcal U_2$ we have
\begin{equation}\label{e76}
\begin{split}a(x_0)
C_{N,s}\int_{U} \frac{\psi(x_0)-\psi(y)}{|x_0-y|^{N+2s}}d y
\\=a(x_0) C_{N,s}\left\{\int_{U}
\frac{\varphi(x_0)-\varphi(y)}{|x_0-y|^{N+2s}}d y +
\int_{\re^N\setminus U} \frac{ h(x_0)- h(y)}{|x_0-y|^{N+2s}}d
y\right\} \\= a(x_0) C_{N,s}\Big\{\int_{\mathcal{U}_1}
\frac{\varphi(x_0)-\varphi(y)}{|x_0-y|^{N+2s}}d y +
\int_{\re^N\setminus\mathcal{U}_1} \frac{ h(x_0)-
h(y)}{|x_0-y|^{N+2s}}d y\\ - \int_{\mathcal{U}_2}
\frac{h(x_0)-h(y)}{|x_0-y|^{N+2s}}d y + \int_{\mathcal{U}_2}
\frac{\varphi(x_0)- \varphi(y)}{|x_0-y|^{N+2s}}dy\Big\}\geq \tilde
C\,.
\end{split}
\end{equation}

\smallskip

Now, take any $x_0 \in \re^N$ with $|x_0|<\bar R$, any open
bounded subset $U\subset \re^N$ with $x_0\in U$, and any test
function $\varphi\in C^2(\re^N)$ such that
$h(x_0)-\varphi(x_0)\leq h(x)-\varphi(x)$ for all $x\in U$. Let
$\psi$ be defined by \eqref{e73}\,. Note that \eqref{e67} gives
\begin{equation}\label{e77} h(x_0)=
W(x_0)\,.
\end{equation}
For any $0< \epsilon< \bar R- R_0$ we have $\mathcal U_1:= U\cap
B_{\bar R-\epsilon}\subset \overline{B}_{\bar R}, x_0\in \mathcal
U_1$. Moreover,
\begin{equation}\label{e75}\varphi(x)\leq W(x) \quad \textrm{for all}\;\; x\in
\mathcal U_1, \quad \varphi(x_0)=W(x_0)\,.\end{equation} So, from
\eqref{e69} with $\Omega'=\mathcal U_1$ we get
\begin{equation}\label{e78}
a(x_0) C_{N,s}\left\{\int_{\mathcal U_1}
\frac{\varphi(x_0)-\varphi(y)}{|x_0-y|^{N+2s}}d y +
\int_{\re^N\setminus\mathcal U_1}
\frac{W(x_0)-W(y)}{|x_0-y|^{N+2s}}d y\right\}\geq  \mu_0\,.
\end{equation}
Due to \eqref{e77} and \eqref{e78}, since $h\leq W$ in $\re^N$, we
have
\begin{equation}\label{e79}
a(x_0) C_{N,s}\left\{\int_{\mathcal U_1}
\frac{\varphi(x_0)-\varphi(y)}{|x_0-y|^{N+2s}}d y +
\int_{\re^N\setminus \mathcal U_1} \frac{h(x_0)-
h(y)}{|x_0-y|^{N+2s}}d y\right\}\geq  \mu_0\,.
\end{equation}
Set $\mathcal U_2:= U\cap \big(\re^N\setminus B_{\bar
R-\epsilon}\big)$.\, In view of \eqref{e79}, since
$\varphi(x_0)-\varphi(y)\geq h(x_0)-h(y)$ for all $y\in \mathcal
U_2$ we have
\begin{equation}\label{e80}
\begin{split}a(x_0)
C_{N,s}\int_{U} \frac{\psi(x_0)-\psi(y)}{|x_0-y|^{N+2s}}d y
\\=a(x_0) C_{N,s}\left\{\int_{U}
\frac{\varphi(x_0)-\varphi(y)}{|x_0-y|^{N+2s}}d y +
\int_{\re^N\setminus U} \frac{ h(x_0)- h(y)}{|x_0-y|^{N+2s}}d
y\right\} \\= a(x_0) C_{N,s}\Big\{\int_{\mathcal{U}_1}
\frac{\varphi(x_0)-\varphi(y)}{|x_0-y|^{N+2s}}d y +
\int_{\re^N\setminus\mathcal{U}_1} \frac{ h(x_0)-
h(y)}{|x_0-y|^{N+2s}}d y\\ - \int_{\mathcal{U}_2}
\frac{h(x_0)-h(y)}{|x_0-y|^{N+2s}}d y + \int_{\mathcal{U}_2}
\frac{\varphi(x_0)- \varphi(y)}{|x_0-y|^{N+2s}}dy\Big\}\geq
\mu_0\,.
\end{split}
\end{equation}
From \eqref{e76} and \eqref{e80} the claim follows. Therefore,
\[h:= \bar C h\quad \textrm{in}\;\; \re^N\,,\]
with $\bar C\geq \max\left\{\frac 1{\mu_0}, \frac 1{\tilde
C}\right\}\,,$ is a supersolution of equation \eqref{e63};
moreover, it is immediately seen that it satisfies \eqref{e64}.
\hfill $\square$

\section{Proofs of existence and uniqueness results} \label{pp}
\setcounter{equation}{0}

To begin with, let us show the next quite standard comparison
principle.
\begin{proposition}\label{prop1a} Let assumptions $(H_0), (H_1)$
be satisfied. Let $u$ be a subsolution of problem \eqref{e1}, let
$v$ be a supersolution of problem \eqref{e1}. Suppose that both
\[\limsup_{|x|\to \infty}(u-v)\,\leq \,0 \quad \textrm{uniformly for}\,\, t\in [0,T]\,. \]
Then
\[u\leq v\quad \textrm{in}\;\; S_T\,.\]
\end{proposition}

\noindent{\it Proof\,.} Set $w:= u- v$. Let $\epsilon>0$. Then
there exists $R_\epsilon>0$ such that
\[
|w(x,t)|\leq \epsilon \quad \textrm{for all}\;\; x\in
\re^N\setminus B_{R_\epsilon}, t\in [0, T]\,.
\]
Hence $w$ is a subsolution of problem
\begin{equation}\label{e56}
\left\{
\begin{array}{ll}
\,\pa_t v = -a\, (-\Delta)^{s}v  + c v
&\textrm{in}\,\,B_{R_\epsilon}\times (0, T]
\\& \\
\textrm{ }v \, = \epsilon & \textrm{in\ \ } \big(\re^N\setminus
B_{R_\epsilon} \big)\times (0, T] \,
\\& \\
\textrm{ }v \, = 0 & \textrm{in\ \ } \re^N\times \{0\} \,.
\end{array}
\right.
\end{equation}
Moreover, it is easily seen that the function
$$z(x,t):= \epsilon \, e^{\|c\|_\infty t}\quad \big(x\in \re^N, t\in [0, T]\big) $$
is a supersolution of problem \eqref{e56}. By the comparison
principle (see Remark \ref{oss1}),
\begin{equation}\label{e57}
w\leq  z \quad \textrm{in}\;\; \re^N\times [0, T]\,.
\end{equation}
Similarly, it can be shown that
\begin{equation}\label{e58}
w\geq - z \quad \textrm{in}\;\; \re^N\times [0, T]\,.
\end{equation}
Letting $\epsilon \to 0^+$, we get $w=0$ in $\re^N\times [0,
T]\,.$ Hence the proof is complete. \hfill $\square$

\medskip
\medskip

Let us prove Theorem \ref{thm1}. Hereafter, $\{\zeta_j\}\subset
C^\infty_c(B_j)$ will be a sequence of functions such that
\begin{equation}\label{e42}
0\leq \zeta_j\leq 1, \;\; \zeta_j\equiv 1\,\,\, \textrm{in}\;\;\,
B_{j/2}\,\;\, \textrm{for each}\,\, j\in \ene\,.
\end{equation}

\smallskip

\noindent{\it Proof of Theorem \ref{thm1}\,.\,\,} For any $j\in
\ene$ let $u_j$ be the unique solution (see Remark \ref{oss1}) of
the problem
\begin{equation}\label{e12} \left\{
\begin{array}{ll}
\,\pa_t u = -a\, (-\Delta)^{s}u  + c u + f &\textrm{in}\,\,
B_j\times (0, T]
\\& \\
\textrm{ }u \, = g & \textrm{in\ \ } \big(\re^N\setminus
B_j\big)\times (0, T] \,
\\& \\
\textrm{ }u \, = u_{0, j}& \textrm{in\ \ } \re^N\times \{0\} \,,
\end{array}
\right.
\end{equation}
where
\[u_{0,j}(x):= \zeta_j(x) u_0(x) + [1-\zeta_j(x)]g(0)\quad \textrm{for all}\;\; x\in \overline{B}_j\,.\]

It is easily seen that the function
\[\overline{v}(x,t):= C e^{\b t} \quad \big((x,t)\in \re^N\times [0,T]\big)\]
is a supersolution of problem \eqref{e12} for any $j\in \ene$,
provided that
\[\b\geq
1+\|c \|_\infty,\;\; C \geq \max\{\|f\|_\infty, \|g\|_\infty, \|u_0\|_\infty\}\,.\]
Thus, by the comparison principle (see Remark \ref{oss1}),
\begin{equation}\label{e13}
u_j(x,t)\leq \bar v(x,t) \quad \textrm{for all}\;\; (x,t)\in
\re^N\times [0,T]\,.
\end{equation}
Furthermore, the function
\[\underline{v}(x,t):= - C e^{\b t} \quad \big((x,t)\in \re^N\times [0,T]\big)\]
is a subsolution of problem \eqref{e12} for any $j\in \ene$. Thus,
by the comparison principle,
\begin{equation}\label{e14}
u_j(x,t)\geq \underline{v}(x,t) \quad \textrm{for all}\;\;
(x,t)\in \re^N\times [0,T]\,.
\end{equation}
From \eqref{e13}-\eqref{e14} we obtain
\begin{equation}\label{e15}
\big|u_j(x,t)\big|\leq C e^{\b  T}=: K_T \quad \textrm{for all}\;\;
(x,t)\in \re^N\times [0,T]\,.
\end{equation}
By the a priori estimates recalled in Remark \ref{oss2}-$(iii)$
and usual compactness arguments, there exists a subsequence
$\{u_{j_k}\}\subset \{u_j\}$ and a function $u\in C(S_T)$ such
that
\[u:= \lim_{k\to \infty} u_{j_k}\quad \textrm{uniformly in}\;\; D\times [\tau, T]\,,\]
for any compact subset $D\subset \re^N$ and for any $\tau\in (0,
T)$. For simplicity we still denote $\{u_{j_k}\}$ by $\{u_j\}\,.$
In view of stability properties of viscosity solutions under local
uniform convergence, the function $u$ is a solution of equation
\[\pa_t u = -a\, (-\Delta)^{s}u  + c u + f \quad \textrm{in}\,\,
\re^N\times (0, T]\,.\]

{\bf Claim 1:} We have that
\[\lim_{t\to 0^+} u(x,t)\,=\, u_0(x) \quad \textrm{for any}\;\; x\in \re^N\,.\]

In fact, let $x_0\in \re^N$. Take $j_0\in \ene$ so large that
$x_0\in B_{j_0/2}.$ In view of the definition of $\{\zeta_j\}$
(see \eqref{e42}) there exists $\d_0>0$ such that for any $j\geq
j_0$
\begin{equation}\label{e43}
u_j(x,0)=u_{0,j}(x)=u_0(x)\quad \textrm{for all}\;\; x\in
B_{\d_0}(x_0)\,.
\end{equation}
Since $u_0\in C(\re^N)$, for any $0<\epsilon<1$ there exists $\d\in
(0, \d_0)$ such that
\begin{equation}\label{e44}
-\epsilon < u_0(x) - u_0(x_0)<\epsilon \quad \textrm{for all}\;\;
x\in B_{\d}(x_0)\,.
\end{equation}
From \eqref{e43}, \eqref{e44} it follows that for any $0<\epsilon<1$
and any $j\geq j_0$ there holds
\begin{equation}\label{e45}
-\epsilon <  u_j(x, 0)- u_0(x_0) <\epsilon \quad \textrm{for all}\;\;
x\in B_{\d}(x_0)\,.
\end{equation}

Let
\[\chi(x;x_0)\equiv \chi(x):=|x-x_0|^2\quad \big(x\in \re^N\big)\,.\]
Since $\chi \in C^\infty(\re^N)$, we have that $(-\Delta)^s
\chi\in C(\re^N)$\,. Define
\[h(x,t):= \big[\chi(x) + A t \big]e^{\eta t}\,\quad (x\in \re^N, t\in [0, \d]),\]
\[\overline v(x,t):=M h(x,t) + u_0(x_0) + \epsilon\, \quad (x\in \re^N, t\in [0, \d])\,,\]
where $A>0, \eta>0, M$ are constants to be determined. We have
that
\[\pa_t \overline v(x,t) =  M[Ae^{\eta t} + \eta h(x,t)]\quad (x\in \re^N, t\in [0,\d]),\]
whereas
\[-a(x)(-\Delta)^s \overline v(x,t) + c(x) \overline v(x,t) + f(x)\]\[
\leq  M e^{\eta \d}\max_{\overline B_\d(x_0)}\big| a (-\Delta)^s
\chi \big| +  M c(x) h(x,t) +\|c\|_\infty(\|u_0\|_\infty
+1)+\|f\|_\infty\quad \big(x\in \re^N, t\in [0,\d]\big).
\]
Therefore,
\begin{equation}\label{e46}
\pa_t \overline v(x,t) \geq -a(x)(-\Delta)^s \overline v(x,t) +
c(x) \overline v(x,t) + f(x) \quad \big(x\in \re^N, t\in
[0,\d]\big),
\end{equation}
if
\begin{equation}\label{e47}
\eta\geq \|c\|_\infty, \quad A\geq \max_{\overline
B_\d(x_0)}\big|a (-\Delta)^s \chi
\big|+\|c\|_\infty(\|u_0\|_\infty+1)+\|f\|_\infty\,.
\end{equation}

Furthermore, since
\[h(x,t)\geq \d^2\quad \textrm{for all}\;\; x\in \re^N\setminus B_\d(x_0), t\in [0,\d],\]
it easily follows that
\begin{equation}\label{e48}
\overline v(x,t) \geq u_j(x,t) \quad \textrm{for all}\;\; x\in
\re^N\setminus B_\d(x_0), t\in [0,\d]\,,
\end{equation}
if
\begin{equation}\label{e49}
M\geq \frac{2K_T}{\d^2}\,.
\end{equation}

From \eqref{e45} we get
\begin{equation}\label{e50}
\overline v(x,0)\geq u_j(x,0)\quad \textrm{for all}\;\; x\in
B_\d(x_0)\,,
\end{equation}
while
\begin{equation}\label{e51}
\overline v(x,0)\geq  M\d^2+ u_0(x_0)\geq u_j(x,0)\quad
\textrm{for all}\;\; x\in \re^N\setminus B_\d(x_0),
\end{equation}
due to \eqref{e49}.

Suppose that \eqref{e47}, \eqref{e49} hold. Then, by \eqref{e46},
\eqref{e48}, \eqref{e50}, \eqref{e51}, for any $j\in \ene, j>R$
the function $\overline v$ is a supersolution (in the sense of
Definition \ref{defsolp}) of problem
\begin{equation}\label{e52} \left\{
\begin{array}{ll}
\,\pa_t v = -a\, (-\Delta)^{s}v  + c v + f
&\textrm{in}\,\,B_\d(x_0)\times (0,\d]
\\& \\
\textrm{ }v \, = u_j & \textrm{in\ \ } \big(\re^N\setminus
B_\d(x_0) \big)\times (0, \d] \,
\\& \\
\textrm{ }v \, = u_j & \textrm{in\ \ } \re^N\times \{0\} \,,
\end{array}
\right.
\end{equation}
while $u_j$ is a solution of the same problem. By the comparison
principle (see Remark \ref{oss1}) we obtain
\begin{equation}\label{e53}
u_j \leq \overline v \quad \textrm{in}\;\; B_\d(x_0)\times
(0,\d]\,.
\end{equation}

Define
\[\underline v(x,t):= - M h(x,t) + u_0(x_0) - \epsilon \quad (x\in \re^N, t\in [0, \d])\,;\]
suppose that \eqref{e47} and \eqref{e49} hold. By the same
arguments as above, we can show that there holds
\begin{equation}\label{e53a}
u_j \geq \underline v \quad \textrm{in}\;\; B_\d(x_0)\times
(0,\d]\,.
\end{equation}
Inequalities \eqref{e53}-\eqref{e53a} yield
\begin{equation}\label{e54}
-M h(x, t) - \epsilon\leq u_j(x,t) - u_0(x_0) \,\leq\,  M h(x,
t)+\epsilon\,
\end{equation}
for all $x\in B_\d(x_0), t\in [0,\d]\,.$ Letting $j\to \infty$,
thus we obtain
\begin{equation}\label{e55}
-M h(x, t) - \epsilon\leq u(x,t) - u_0(x_0) \,\leq\, M h(x,
t)+\epsilon\,
\end{equation}
for all $x\in B_\d(x_0), t\in (0,\d]\,.$ Letting $x\to x_0, t\to
0^+$, and then $\epsilon \to 0^+$, we get that $\lim_{x\to x_0}
u(x,t)=u_0(x_0)$. Hence the Claim 1 has been shown.

\bigskip

\noindent {\bf Claim 2:} We have that
\[ \lim_{|x|\to \infty} u(x,t)\,=\, g(t) \quad \textrm{uniformly with respect to}\;\; t\in [0,T]\,.\]

In fact, fix any $\t_0\in [0, T], 0<\epsilon<1$. Since $g\in C([0,
T])$, there exists $\d\in (0,1)$ such that
\begin{equation}\label{e25}
g(t_0)-\epsilon \leq g(t)\leq g(t_0) +\epsilon \quad \textrm{for
any}\;\; t\in [\underline t_\d, \overline t_\d],
\end{equation}
where
\[\underline t_\d:= \max\{t_0-\d, 0\}, \quad \overline t_\d:= \min\{t_0+\d, T\}\,.\]
Clearly, $\d=\d(\epsilon)$ does not depend on $t_0$. Furthermore,
due to \eqref{e8}, there exists $R_\epsilon>0$ such that
\begin{equation}\label{e31a}
g(0)-\epsilon \leq u_0(x) \leq g(0) +\epsilon \quad \textrm{for
all}\;\; x\in \re^N\setminus B_{R_\epsilon}\,.
\end{equation}
Let $R\geq \max\{R_0, R_\epsilon\}$ with $R_0$ given by
Proposition \ref{prop1}; set
\[N_j^R:= B_j\setminus B_{R} \quad \textrm{for any}\;\; j>R\,.\]

Define
\begin{equation}\label{e26}
\underline w(x,t):= - M V(x) e^{\eta t} -\l (t-t_0)^2 + g(t_0)
-\epsilon\quad \textrm{for all}\;\; (x,t)\in \re^N\times [0, T]\,,
\end{equation}
where $M>0, \eta>0, \l>0$ are constants to be chosen in the
sequel, while  $V(x)\equiv V(|x|)$ is the supersolution given by
Proposition \ref{prop1}\,.

\medskip

In view of Proposition \ref{prop1}, we have
\[-a(x)(-\Delta)^s
\underline w + c(x) \underline w \geq M e^{\eta t} - M \, c(x)
V(x)\,e^{\eta t} - \|c\|_\infty(\|g\|_\infty +\l +1)\quad
\textrm{for all}\;\; x\in N_j^R, t\in [0, T]\,.
\]
Therefore,
\begin{equation}\label{e37a}
\begin{split}
\pa_t \underline w  + a (-\Delta)^s\underline w - c \underline w
-  f \hspace{3 cm}\\
\leq - \eta  M V e^{\eta t} - 2 \l(t-t_0) - M e^{\eta
t} + c M V e^{\eta t} \hspace{1.8 cm} \\
+ \|c\|_\infty\big(\|g\|_\infty + \l +1\big)+\|f\|_\infty \leq 0
\quad \textrm{in}\;\; N_j^R\times (0, T]\,, \hspace{1 cm}
\end{split}
\end{equation}
if we take
\begin{equation}\label{e27}
\eta\geq \|c \|_\infty\,,
\end{equation}
\begin{equation}\label{e28}
M \geq 2 \l + \|f\|_\infty + \|c\|_\infty(\|g\|_\infty+
\l+1)+\|f\|_\infty\,.
\end{equation}
In view of \eqref{e15}, we obtain
\begin{equation}\label{e29}
\underline w(x,t) \leq - M V(R) + \|g\|_\infty \leq -K_T\leq
u_j(x,t) \quad \textrm{for all}\;\; x\in \overline{B}_R, t\in
(\underline t_\d, \overline t_\d)\,,
\end{equation}
if
\begin{equation}\label{e36}
M\geq \frac{\|g\|_\infty + K_T}{V(R)}\,.
\end{equation}

From \eqref{e25} we have
\begin{equation}\label{e30}
\underline w(x,t) \leq g(t) \quad \textrm{for all}\;\; x\in
\re^N\setminus B_j, t\in (\underline t_\d, \overline t_\d)\,.
\end{equation}

Suppose that $\underline t_\d=0$ (note that this is always the case when
$t_0=0$). From \eqref{e31a} and \eqref{e25}
\begin{equation}\label{e31}
\underline w(x,0) \leq g(0)-\epsilon\leq u_j(x,0)=u_{0,j}(x) \quad
\textrm{for all}\;\; x\in \re^N\setminus B_{R}\,;
\end{equation}
while
\begin{equation}\label{e32}
\underline w(x,0)\leq -M V(R_\epsilon) +\|g\|_\infty \leq -K_T\leq
u_j(x,0)\quad \textrm{for all}\;\; x\in B_{R}\,,
\end{equation}
provided that \eqref{e36} holds.

Suppose that $\underline t_\d>0$\,. It follows from \eqref{e15}
that
\begin{equation}\label{e34}
\underline w(x,\underline t_\d) \leq - \l \d^2 +\|g\|_\infty \leq
-K_T\leq u_j(x,\underline t_\d)\quad \textrm{for all}\;\; x\in
\re^N, t\in (\underline t_\d, \overline t_\d)\,,
\end{equation}
if
\begin{equation}\label{e35}
\l\geq \frac{\|g\|_\infty + K_T}{\d^2}\,.
\end{equation}

Now, suppose that \eqref{e27}, \eqref{e28}, \eqref{e36},
\eqref{e35} hold. By \eqref{e37a}, \eqref{e29}, \eqref{e30},
\eqref{e31}, \eqref{e32}, \eqref{e34}, for any $j\in \ene, j>R$,
the function $\underline w$ is a subsolution (in the sense of
Definition \ref{defsolp}) of problem
\begin{equation}\label{e37} \left\{
\begin{array}{ll}
\,\pa_t v = -a\, (-\Delta)^{s}v  + c v + f
&\textrm{in}\,\,N^R_j\times (\underline t_\d, \overline t_\d]
\\& \\
\textrm{ }v \, = u_j & \textrm{in\ \ } \big(\re^N\setminus N^R_j
\big)\times (0, T] \,
\\& \\
\textrm{ }v \, = u_j & \textrm{in\ \ } \re^N\times \{0\} \,,
\end{array}
\right.
\end{equation}
while $u_j$ is a solution of the same problem. By the comparison
principle (see Remark \ref{oss1}) we obtain
\begin{equation}\label{e38}
\underline w \leq u_j \quad \textrm{in}\;\; N^R_j\times
(\underline t_\d, \overline t_\d]\,.
\end{equation}

Define
\begin{equation}\label{e26a}
\overline w(x,t):= M V(x)
e^{\eta t} + \l (t-t_0)^2 + g(t_0)+\epsilon\quad \textrm{for
all}\;\; x\in \re^N, t\in [0,T]\,;
\end{equation}
suppose that
\eqref{e27}, \eqref{e28}, \eqref{e36}, \eqref{e35}\,. By the same
arguments as above, we can show that there holds
\begin{equation}\label{e38a}
\overline w \geq u_j \quad \textrm{in}\;\; N^R_j\times (\underline
t_\d, \overline t_\d]\,.
\end{equation}
From \eqref{e38} and \eqref{e38a} we get
\begin{equation}\label{e39}
-M V(x) e^{\eta t} - \l (t-t_0)^2-\epsilon\leq u_j(x,t)
- g(t_0) \,\leq\, M V(x) e^{\eta t} + \l (t-t_0)^2 +\epsilon\,
\end{equation}
for all $x\in N_j^R, t\in (\underline t_\d, \overline t_\d]\,.$
Choosing $t=t_0$ in \eqref{e39} and letting $j\to \infty$, we
obtain
\begin{equation}\label{e40}
- M V(x) e^{\eta T} -\epsilon\leq u(x,t_0) - g(t_0) \,\leq\, M
V(x) e^{\eta T} +\epsilon\,\quad \textrm{for all}\;\; x\in
\re^N\setminus B_R\,.
\end{equation}
From \eqref{e40} it follows that
\begin{equation}\label{e41}
\sup_{t_0\in [0, T]}\big|u(x,t_0) - g(t_0) \big| \leq \overline C
V(x) + \epsilon \quad \textrm{for all}\;\; x\in \re^N\setminus
B_R\,,
\end{equation}
where $\overline C:= M e^{\eta T}$. Due to \eqref{e41} and
\eqref{e20}, letting $|x|\to \infty, \epsilon\to 0^+$, we obtain
\eqref{e9}. Hence the Claim 2 has been shown.

Finally, this solution is unique, due to Proposition
\ref{prop1a}\,. \hfill $\square$

\medskip

\bigskip
\medskip

Now we prove Theorem \ref{thm1a}. We follow the same line of
arguments of the proof of Theorem \ref{thm1a}, but there is an
important difference. In fact, we need to substitute the estimate
\eqref{e15}, which is dependent on $T$, by another one independent
of $T$. In order to obtain such better estimate we use the
supersolution $h$ constructed in Proposition \ref{prop2}.

\medskip

\noindent {\it Proof of Theorem \ref{thm1a}\,.} Arguing as in the
proof of Theorem \ref{thm1} we construct the sequence $u_j(x,t)$
of solutions of problem \eqref{e12} with $T=\infty$. Let $h(x)$ be
the supersolution provided by Proposition \ref{prop2}. Then
obviously
\begin{equation}\label{e110}
V_0(x):= h(x) - \inf_{\re^N} h + 1
\end{equation}
is also a supersolution of \eqref{e63} and $V_0(x)\geq 1$. Let
$B:=\max\{\|f\|_\infty, \| u_0\|_\infty, \|g\|_\infty\}$. Since $c\leq 0,$ we have that $B V_0$ is
a supersolution of problem \eqref{e12}, while $-B V_0$ is a
subsolution of \eqref{e12}. Thus, by the comparison principle,
\begin{equation}\label{e111}
|u_j|\, \leq \, B V_0 \quad \textrm{in}\;\;\, B_j\times
(0,\infty)\,.
\end{equation}
Passing to the limit as $j\to\infty$ we obtain that
\begin{equation}\label{e112}
|u|\leq B V_0\leq \check C:= B\| V_0\|_\infty \quad
\textrm{in}\;\; \re^N\times (0,\infty)\,.
\end{equation}
Note that estimate \eqref{e112} substitutes estimate \eqref{e15}
which is depending on $T$. Now, consider the functions $\underline
w$ and $\overline w$ defined in \eqref{e26} and in \eqref{e26a},
respectively. Suppose that $\eta=0$,
$$\l\geq \frac{\|g\|_\infty + \check C}{\d^2},\, M\geq \frac{\|g\|_\infty + \check C}{V(R)},$$
and \eqref{e28} holds. Note that $M$ and $\l$ do not depend on
$T$. Since $c\leq 0$ and \eqref{e112} holds, by the same arguments
as in the proof of Theorem \ref{thm1} we can infer that for any
$\epsilon>0$
\begin{equation}\label{e41c}
\sup_{t_0\in [0, \infty)}\big|u(x,t_0) - g(t_0) \big| \leq
 M V(x) + \epsilon \quad \textrm{for all}\;\; x\in
\re^N\setminus B_R\,.
\end{equation}
Thanks to \eqref{e41c} and \eqref{e20}, letting $|x|\to \infty,
\epsilon\to 0^+$, we obtain \eqref{e9a}. Finally, this solution is
unique, due to Proposition \ref{prop1a}\,. This completes the
proof\,. \hfill $\square$

\bigskip
\medskip

We have the next quite standard comparison principle.
\begin{proposition}\label{prop7}
Let assumptions $(H_0), (H_1)$ be satisfied. Suppose that $c\leq
0$ in $\re^N$. Let $u$ be a subsolution and $v$ a supersolution to
equation \eqref{e3} such that
\[\limsup_{|x|\to\infty} (u-v)\leq 0\,.\] Then $$u\leq v \quad \textrm{in}\;\; \re^N\,.$$
\end{proposition}

\noindent{\it Proof\,.} Set $w:= u_1-u_2$. Let $\epsilon>0$. Then
there exists $R_\epsilon>0$ such that
\[
|w(x)|\leq \epsilon \quad \textrm{for all}\;\; x\in \re^N\setminus
B_{R_\epsilon}\,.
\]
Hence $w$ is a subsolution of problem
\begin{equation}\label{e85}
\left\{
\begin{array}{ll}
\,-a\, (-\Delta)^{s}v  + c v\,=\,0 &\textrm{in}\,\,B_{R_\epsilon}
\\& \\
\textrm{ }v \, = \epsilon & \textrm{in\ \ } \re^N\setminus
B_{R_\epsilon}\,.
\end{array}
\right.
\end{equation}
Moreover, it is easily seen that the function $z\equiv \epsilon $
is a supersolution of problem \eqref{e85}. So, by the comparison
principle (see Remark \ref{oss2}),
\begin{equation}\label{e86}
w\leq  \epsilon \quad \textrm{in}\;\; \re^N\,.
\end{equation}
Similarly, it can be shown that
\begin{equation}\label{e87}
w\geq - \epsilon \quad \textrm{in}\;\; \re^N\,.
\end{equation}
Letting $\epsilon \to 0^+$, we get $w=0$ in $\re^N\,.$ Hence the
proof is complete. \hfill $\square$

\medskip
\medskip

Now, we prove Theorem \ref{thm2}\,.

\noindent{\it Proof of Theorem \ref{thm2}\,.} Let $\g\in \re$. For
any $j\in \ene$ let $u_j$ be the unique solution (see Remark
\ref{oss2}) of the problem
\begin{equation}\label{e59} \left\{
\begin{array}{ll}
\,a\, (-\Delta)^{s}u  - c u\, =\, f &\textrm{in}\,\, B_j\times (0,
T]
\\& \\
\textrm{ }u \, = \g & \textrm{in\ \ } \re^N\setminus B_j\,\,.
\end{array}
\right.
\end{equation}

We claim that there exists $K>0$ such that for any $j\in \ene$
\begin{equation}\label{e60}
\big|u_j(x)\big|\leq K \quad \textrm{for all}\;\; x\in \re^N\,.
\end{equation}
In fact, let $h=h(x)\equiv h(|x|)$ be the supersolution given by
Proposition \ref{prop2}. Define
$$
\tilde h := C (h +1)\quad \textrm{in}\;\; \re^N\,,
$$
where $C\geq \max\{\g, \|f\|_\infty\}\,.$ It is easily seen that,
for any $j\in \ene,$ $h$ is a supersolution of problem
\eqref{e59}. Therefore, by the comparison principle (see Remark
\ref{oss2}), we get \eqref{e60}, with $K=\|\tilde h\|_\infty\,.$

By the a priori estimates recalled in Remark \ref{oss2}-$(iii)$
and usual compactness arguments, there exists a subsequence
$\{u_{j_k}\}\subset \{u_j\}$ and a function $u\in C(\re^N)$ such
that
\[u:= \lim_{k\to \infty} u_{j_k}\quad \textrm{uniformly in}\;\; D\,,\]
for any compact subset $D\subset \re^N$. For simplicity, we still
denote $\{u_{j_k}\}$ by $\{u_j\}\,.$ In view of stability
properties of viscosity solutions under local uniform convergence,
the function $u$ is a solution of equation
\[a\, (-\Delta)^{s}u - c u \,=\, f \quad \textrm{in}\,\,
\re^N\,.\]

\smallskip

{\bf Claim\,:} The solution $u$ satisfies condition \eqref{e10c}.

In fact, define
\begin{equation}\label{e81}
\underline w(x):= - M h(x) + \g \quad \textrm{for
all}\;\; x\in \re^N\,,
\end{equation}
where $ M>0$ is a constant to be chosen in the sequel\,.

\medskip

In view of Proposition \ref{prop2}, it is easily seen that, if we
take
\[
M \geq \|c\|_\infty\g+\|f\|_\infty\,,
\]
then $\underline w$ is a supersolution of problem \eqref{e59}, for
any $j\in \ene.$ By the comparison principle (see Remark
\ref{oss2}),
\begin{equation}\label{e83}
\underline w \leq u_j\quad \textrm{in}\;\; \re^N\,.
\end{equation}
On the other hand, by the same methods as above, we can show that
\begin{equation}\label{e84}
u_j\leq \overline w \quad \textrm{in}\;\; \re^N\,,
\end{equation}
where
\[\overline w(x):= M h(x) + \g \quad \textrm{for
all}\;\; x\in \re^N\,.\] 

From \eqref{e83}, \eqref{e84} it follows that
\[
- M h + \g \leq u_j\leq  M h + \g\quad
\textrm{in}\;\; \re^N\,.
\]
Letting $j\to \infty$, in view of \eqref{e64} we have that
\eqref{e10} holds. So, the Claim has been shown.

Finally, the uniqueness of the solution $u$ follows from
Proposition \ref{prop7}\,. \hfill $\square$

\section{Asymptotic behaviour of solutions: proofs} \setcounter{equation}{0}
To begin with, we show the next auxiliary result.
\begin{proposition}\label{prop3}
Let assumptions of Theorem \ref{thm1a} be satisfied with $g\equiv
g_1.$ Assume that
\begin{equation}\label{e113}
g_1(t_1)\,\leq\, g_1(t_2)\quad \textrm{for any}\,\, 0\leq
t_1<t_2\,.
\end{equation}
Let $\underline V:= - A V_0$, with
\begin{equation}\label{e114}
A \geq \|g_1\|_\infty + \| f\|_\infty\,
\end{equation}
and $V_0$ defined in \eqref{e110}\,. Let $\underline w$ be the
unique solution, provided by Theorem \ref{thm1a}, of the problem
\begin{equation}\label{e115}
\left\{
\begin{array}{ll}
\,   \, \pa_t u =- a(-\Delta)^s u\, +\,  c u\,+f
&\textrm{in}\,\,\re^N\times (0,\infty)
\\&\\
\textrm{ } u \, = \underline V & \textrm{in\ \ } \re^N\times \{0\}
\end{array}
\right.
\end{equation}
such that
\begin{equation}\label{e116}
\lim_{|x|\to \infty} \underline w(x, t)\,= g_1(t)\quad
\textrm{uniformly for}\;\; t\in [0,\infty)\,.
\end{equation}
Then $t\mapsto \underline w(x,t)$ is nondecreasing, i.e.,
\begin{equation}\label{e117}
\underline w(x, t_1)\leq \underline w(x, t_2)\quad \textrm{for
all}\;\; x\in \re^N, 0 \leq  t_1<t_2\,.
\end{equation}
\end{proposition}

\medskip
\medskip

\noindent{\it Proof of Proposition \ref{prop3}\,.\,\,}  It is
easily seen that $\underline V$ is a subsolution of problem
\eqref{e115}\,. In fact, since $c\leq 0$ and $\underline V<0$, due
to \eqref{e114} we have (in the viscosity sense)
\[-a(-\Delta)^s \underline V + c \underline V + f \geq  A - \|f\|_\infty\geq 0=\pa_t \underline V \quad \textrm{in}\,\, \re^N\times (0,\infty)\,.\]
Moreover,
\[\underline V -\underline w =0\quad \textrm{in}\,\, \re^N\times \{0\},\]
and by \eqref{e114} and \eqref{e116},
\[\limsup_{|x|\to \infty} [\underline V(x) -  \underline w(x, t)]\, \leq 0\quad \textrm{uniformly for}\;\; t\in [0,\infty)\,.\]
Since $\underline w$ is a solution of problem \eqref{e115}, by Proposition \ref{prop1a},
\begin{equation}\label{e118}
\underline V(x)=\underline w(x,0) \, \leq\, \underline w(x,t)
\quad \textrm{for all}\;\; x\in \re^N, t>0\,.
\end{equation}

In order to show \eqref{e117}, take any $t_0>0$ and define
\[\tilde w(x,t):= \underline w (x, t+ t_0)\quad \textrm{for all}\;\; x\in \re^N, t>0\,.\]
Note that both $\underline w$ and $\tilde w$ satisfy the equation
\[\pa_t v - a(-\Delta)^s v - c v \,=\, f\quad \textrm{in}\;\; \re^N\times (0,\infty)\,. \]
Moreover, from \eqref{e118} we obtain that
\begin{equation}\label{e119}
\tilde w(x,0) \geq \underline w (x,0)\quad \textrm{for all}\;\;
x\in \re^N\,.
\end{equation}
In addition, due to \eqref{e113},
\[\lim_{|x|\to \infty}[ \tilde w(x,t) - \underline w(x,t) ] = \tilde g_1(t+t_0)-g_1(t)\geq 0\quad \textrm{uniformly for}\,\, t\in [0,\infty)\,. \]
Thus, by Proposition \ref{prop1a},
\[  \tilde w(x,t)\geq \underline w(x,t)\quad \textrm{for all}\;\; x\in \re^N, t>0\,.\]
Hence the conclusion follows. \hfill $\square$

\medskip
\medskip
Similarly, we can show the next result.
\begin{proposition}\label{prop4}
Let assumptions of Theorem \ref{thm1a} be satisfied with $g\equiv
g_2\,.$ Assume that
\begin{equation}\label{e120}
g_2(t_1)\,\geq\, g_2(t_2)\quad \textrm{for any}\,\, 0\leq
t_1<t_2\,.
\end{equation}
Let $\overline V:= - A V_0$, where $V_0$ is defined in
\eqref{e121} and $A$ in \eqref{e114}\,.

Let $\overline w$ be the unique solution, provided by Theorem
\ref{thm1a}, of the problem
\begin{equation}\label{e121}
\left\{
\begin{array}{ll}
\,   \, \pa_t u =- a(-\Delta)^s u\, +\,  c u\,+f
&\textrm{in}\,\,\re^N\times (0,\infty)
\\&\\
\textrm{ } u \, = \overline V & \textrm{in\ \ } \re^N\times \{0\}
\end{array}
\right.
\end{equation}
such that
\begin{equation}\label{e122}
\lim_{|x|\to \infty} \overline w(x, t)\,= g_2(t)\quad
\textrm{uniformly for}\;\; t\in [0,\infty)\,.
\end{equation}
Then $t\mapsto \overline w(x,t)$ is nonincreasing, i.e.,
\begin{equation}\label{e123}
\overline w(x, t_1)\geq \overline w(x, t_2)\quad \textrm{for
all}\;\; x\in \re^N, 0 \leq  t_1<t_2\,.
\end{equation}
\end{proposition}

\medskip
\medskip

Now we prove the next result.
\begin{proposition}\label{prop5}
Let assumptions of Theorem \ref{thm1a} be satisfied. Let $g_1\in
C([0,\infty))\cap L^\infty((0,\infty))$ with
\begin{equation}\label{e126}
g_1(t)\leq g(t)\quad \textrm{for all}\;\, t\in [0, \infty),
\end{equation}
\begin{equation}\label{e133}
\lim_{t\to \infty} g_1(t) = \lim_{t\to \infty} g(t)\,;
\end{equation}
suppose that \eqref{e113} is satisfied. Let $\underline w$ be
given by Proposition \ref{prop3}, also supposing that
\begin{equation}\label{e124}
A\geq \| u_0\|_\infty\,.
\end{equation}
Then
\begin{equation}\label{e125}
\underline w(x,t) \leq u(x,t) \quad \textrm{for all}\,\, x\in
\re^N, t>0\,.
\end{equation}
\end{proposition}

\noindent{\it Proof\,.\,\,} Let $z:=\underline w - u$. Note that
$z$ solves equation
\[\pa_t z  = -a (-\Delta)^s z + c z \quad \textrm{in}\;\; \re^N\times (0,\infty)\,.\]
In view of \eqref{e124} we have
\[ z(x,0)= \underline V(x) - u_0(x) \leq 0\quad \textrm{for all}\;\; x\in \re^N\,.\]
Moreover, from \eqref{e126} we obtain
\[ \lim_{|x|\to \infty} z(x, t) = g_1(t)- g(t)\leq 0 \quad \textrm{uniformly for}\;\; t\in [0,\infty)\,. \]
Hence, by Proposition \ref{prop1a},
\[ z\leq 0\quad \textrm{for all}\;\; x\in \re^N, t>0\,.\]
This completes the proof. \hfill $\square$

\bigskip
Analogously to Proposition \ref{prop5}, the next result can be
shown.

\begin{proposition}\label{prop6}
Let assumptions of Theorem \ref{thm1a} be satisfied. Let $g_2\in
C([0,\infty))\cap L^\infty((0,\infty))$ with
\begin{equation}\label{e130}
g_2(t)\geq g(t)\quad \textrm{for all}\;\, t\in [0, \infty),
\end{equation}
\begin{equation}\label{e134}
\lim_{t\to \infty} g_2(t) = \lim_{t\to \infty} g(t)\,;
\end{equation}
suppose that \eqref{e120} is satisfied. Let $\overline w$ be given
by Proposition \ref{prop4}, also supposing that \eqref{e124}
holds. Then
\begin{equation}\label{e131}
\overline w(x,t) \geq u(x,t) \quad \textrm{for all}\,\, x\in
\re^N, t>0\,.
\end{equation}
\end{proposition}

\bigskip
\bigskip

Now we are in position to prove Theorem \ref{thm3}\,.

\smallskip

\noindent{\it Proof of Theorem \ref{thm3}\,.\,\,} Keep the same
notation as in Propositions \ref{prop3}-\ref{prop6}. In view of
\eqref{e117} and \eqref{e123}, we can define
\begin{equation}\label{e136} \underline W(x):= \lim_{t\to \infty}
\underline w(x,t), \quad \overline W(x):= \lim_{t\to \infty}
\overline w(x,t)\quad \textrm{for any}\;\; x\in
\re^N\,.\end{equation} Observe that the constant $C$ in Remark
\ref{oss1} do not depend on $T$, since $a, c, f$ does not depend
on $t$. Consequently we have that $\underline w\to \underline W,\,
\overline w\to \overline W$ as $t\to \infty$ uniformly in each
compact subset of $\re^N$; thus, $\underline W, \overline W\in
C(\re^N)$\,. We claim that both $\underline W$ and $\overline W$
solve
\begin{equation}\label{e135}
a(-\Delta)^s u - c u \,=\, f\quad \textrm{in}\;\;
\re^N\,.
\end{equation}
In fact, we limit ourselves to show that $\underline W$ is a
subsolution of equation \eqref{e135}, since the remaining part of
the claim follows analogously.

Now, let $\{t_n\}\subset (0,\infty)$ be a sequence with $t_n\to
\infty$ as $n\to \infty$\,. Set
\[\underline w_n(x):= \underline w(x, t_n)\quad (x\in \re^N)\,.\]
Thus, $\underline w_n\to \underline W$ locally uniformly in
$\re^N$ as $n\to \infty.$

Take any bounded subset $U\subset \re^N, x_0\in U,$ take any test
function $\varphi\in C^2(\re^N)$ such that
$$\underline W(x_0) -\varphi(x_0) \geq \underline W(x)-\varphi(x) \quad \textrm{for all}\;\; x\in U\,.$$

\smallskip

Choose $\xi\in C^2(\re^N)$ with \begin{equation}\label{e142} 0\leq
\xi< 1\quad \textrm{if}\,\, x\in \re^N\setminus \{x_0\},\,\,
\xi(x_0)=1\,.\end{equation} Fix any $\epsilon>0.$ So,
\[\underline W(x_0)-[\varphi(x_0)-\epsilon\xi(x_0)] > \underline W(x) -[\varphi(x)-\epsilon \xi(x)]\quad \textrm{for all}\;\; x\in U\setminus\{x_0\}\,.\]
 It is easily seen that there exists
$\bar n=\bar n(\epsilon)\in\ene$ such that for any $n > \bar n$,
for some $x^\epsilon_n\in U$,
\[\underline w_n(x^\epsilon_n) - [\varphi(x^\epsilon_n)-\epsilon \xi(x^\epsilon_n)] \geq \underline w_n(x) - [\varphi(x)-\epsilon \xi(x)]\quad \textrm{for all}\;\; x\in U\,;\]
moreover, for each $\epsilon>0$, $x^\epsilon_n\to x_0$ as $n\to
\infty\,$.

Since $\underline w$ is a solution of \eqref{e115}, due to
Definition \ref{defsoleqp}, we have that
\begin{equation}\label{e141} 0=\partial_t \chi(x^\epsilon_n)\leq -
a(x^\epsilon_n)(-\Delta)^s \chi(x^\epsilon_n) +
c(x^\epsilon_n)\underline w_n(x^\epsilon_n)+
f(x^\epsilon_n),\end{equation} with
\[\chi\equiv \chi_{\epsilon, n}:=\left\{
\begin{array}{ll}
\,\varphi-\epsilon\xi &  \textrm{in}\,\,U
\\& \\
\textrm{ }\underline w_n & \textrm{in\ \ } \re^N\setminus U \,.
\end{array}
\right.
\]
Note that
\begin{equation}\label{e143} (-\Delta)^s \chi(x_n)=
C_{N,s}\left\{\int_{U}\frac{\varphi(x^\epsilon_n)-\epsilon\xi(x^\epsilon_n)-[\varphi(y)-\epsilon\xi(y)]}{|x^\epsilon_n-y|^{N+2s}}
dy + \int_{\re^N\setminus U}\frac{\underline
w_n(x^\epsilon_n)-\underline w(y)}{|x^\epsilon_n-y|^{N+2s}} dy
 \right\}\,. \end{equation}
Since $\varphi, \chi\in C^2(U)$, for any $\epsilon>0 $ we have
\begin{equation}\label{e144}
\begin{split}
\lim_{n\to
\infty}\int_{U}\frac{\varphi(x^\epsilon_n)-\epsilon\xi(x^\epsilon_n)-[\varphi(y)-\epsilon\xi(y)]}{|x^\epsilon_n-y|^{N+2s}}
dy \,=\, \int_{U}\frac{\varphi(x_0)-\varphi(y)}{|x_0-y|^{N+2s}} dy
+ \epsilon\int_{U}\frac{\xi(y)-\xi(x_0)}{|x_0-y|^{N+2s}} dy\,;
\end{split}
\end{equation}
furthermore,
\begin{equation}\label{e145}
\lim_{n\to \infty}\int_{\re^N\setminus U}\frac{\underline
w_n(x^\epsilon_n)-\underline w_n(y)}{|x^\epsilon_n-y|^{N+2s}}
dy\,=\, \int_{\re^N\setminus U}\frac{\underline W(x_0)-\underline
W(y)}{|x_0-y|^{N+2s}} dy\,.\end{equation} From \eqref{e143},
\eqref{e144}, \eqref{e145}, letting $n\to \infty$ in \eqref{e141},
we have, for any $\epsilon>0$,
\[0\leq -a(x_0)(-\Delta)^s \psi(x_0)-a(x_0)\epsilon\int_{U}\frac{\xi(y)-\xi(x_0)}{|x_0-y|^{N+2s}}
dy + c(x_0)\underline w(x_0)+ f(x_0)\,,\] with
\[\psi:=\left\{
\begin{array}{ll}
\,\varphi &  \textrm{in}\,\,U
\\& \\
\textrm{ }\underline W & \textrm{in\ \ } \re^N\setminus U \,.
\end{array}
\right.
\]
Letting $\epsilon\to 0$, the claim follows.

\smallskip

Note that, in view of \eqref{e116}, \eqref{e122}, \eqref{e133},
\eqref{e134}, we can infer that
\[\lim_{|x|\to \infty} \underline W(x)=\overline W(x) =\gamma\,,
\]
where $\g=\lim_{t\to \infty} g(t)\,.$ By Proposition \ref{prop7},
\begin{equation}\label{e151}
\underline W(x)=\overline W(x)\quad \textrm{for all}\;\; x\in
\re^N\,.\end{equation} By \eqref{e125} and \eqref{e131},
\[\underline w(x,t) \leq u(x,t) \leq \overline w(x,t) \quad \textrm{for all}\;\; x\in \re^N, t>0\,.\]
Letting $t\to \infty$, due to \eqref{e136} and \eqref{e151}, we
get the thesis, with $W:= \underline W\equiv\overline W\,.$

 \hfill $\square$

\end{document}